\documentclass[twocolumn]{autart} 

\usepackage{graphicx} 
\usepackage{enumerate}
\usepackage[round, semicolon]{natbib}
\usepackage[letterpaper, left=0.75in, right=0.75in, top=1in, bottom=1in]{geometry}
\usepackage{subcaption}
\usepackage{amsmath} 
\usepackage{amssymb}  
\usepackage{color}
\usepackage[T1]{fontenc}
\usepackage{todonotes}
\usepackage{comment}
\usepackage{wrapfig}
\begin{document}

\begin{frontmatter}

\title{
Decentralized Optimal Coordination of Connected and Automated Vehicles for Multiple Traffic Scenarios\thanksref{footnoteinfo}}

\author[UD]{A M Ishtiaque Mahbub}\ead{mahbub@udel.edu},
\author[UD]{Andreas Malikopoulos}\ead{andreas@udel.edu},
\author[UD]{Liuhui Zhao}\ead{lhzhao@udel.edu}

\address[UD]{Department of Mechanical Engineering, University of Delaware, 126 Spencer Lab, 130 Academy St, Newark DE 19716}

\thanks[footnoteinfo]{This research was supported by ARPAE's NEXTCAR program under the award number DE-AR0000796.}


\begin{abstract}
Connected and automated vehicles (CAVs) provide the most intriguing opportunity to optimize energy consumption and travel time. Several approaches have been proposed in the literature that allow CAVs to coordinate in situations where there is a potential conflict, for example, in signalized intersections, merging at roadways and roundabouts, to reduce energy consumption and optimize traffic flow. In this paper, we consider the problem of coordinating CAVs in a corridor consisting of multiple traffic scenarios. We formulate a two-level optimization problem in which we maximize traffic throughput in the upper-level problem, and derive a closed-form analytical solution that yields the optimal control input for each CAV, in terms of fuel consumption, in the low-level problem. We validate the effectiveness of the solution through simulation under 100\% CAV penetration rate. Fuel consumption and travel time for the vehicles are significantly reduced compared to a baseline scenario consisting of human-driven vehicles. 
\end{abstract}

\end{frontmatter}


\section{Introduction} \label{sec:in}

Connectivity and automation in vehicles provide the most intriguing opportunity for enabling users to better monitor transportation network conditions and make better operating decisions. Several research efforts have been reported in the literature proposing optimization and control approaches for coordinating CAVs at different traffic scenarios that include merging at roadways and roundabouts, crossing intersections, cruising in congested traffic, passing through speed reduction zones, and lane-merging or passing maneuvers. \cite{Dresner2004} presented the use of the reservation scheme to control a single intersection of two roads with vehicles traveling with similar speed on a single direction on each road. Following this effort, similar approaches have been reported in the literature for safe and efficient coordination of CAVs at urban intersections; see \cite{Dresner2008, DeLaFortelle2010, Huang2012, Colombo2014}. \cite{Kim2014} proposed an approach based on model predictive control that allows each vehicle to optimize its movement locally in a distributed manner with respect to any objective of interest. \cite{Colombo2014} constructed the invariant set for the control inputs that ensure lateral collision avoidance. Several papers have also focused on multi-objective optimization problems using a receding horizon control solution either in centralized or decentralized approaches; see \cite{Kamal2013a, Makarem2013, qian2015}.
Lately, a decentralized optimal control framework was presented for coordinating online CAVs in different transportation scenarios, yet without considering state and control constraints; see \cite{Ntousakis:2016aa, Zhao2019CCTAa}, or without considering rear-end collision avoidance constraint; see \cite{Malikopoulos2017}. A detailed discussion of the research efforts in this area can be found in recent survey papers \citep{Malikopoulos2016a, Guanetti2018}.

In this paper, we address the problem of optimally coordinating CAVs in a corridor consisting of multiple traffic scenarios to improve energy consumption and travel time under the hard safety constraints of colision avoidance. We formulate a two-level optimization problem in which we maximize traffic throughput in the upper-level problem, and derive a closed-form analytical solution that yields the optimal control input for each CAV, in terms of fuel consumption, in the low-level problem. In earlier work, we presented a preliminary analysis on coordinating CAVs in a corridor yet without considering the rear-end safety constraint; see \citep{zhao2018} and \cite{Mahbub2019ACC}. Thus, the contribution of this paper is the formulation and analytical solution of an optimal control problem for coordinating CAVs in a corridor consisting of multiple traffic scenarios with the explicit incorporation of the rear-end safety constraint.

The structure of the paper is organized as follows. In Section \ref{sec:pf}, we formulate the problem, provide the modeling framework, and derive the analytical, closed-form solution with rear-end safety constraint. In Section \ref{sec:sim}, we validate the effectiveness of the analytical solution in a simulation environment and conduct a comparison analysis with traditional human-driven vehicles. Finally, we provide concluding remarks and discussion in Section \ref{sec:conc}.

\section{Problem Formulation} \label{sec:pf}
\begin{figure}[!ht]
\centering
\includegraphics[width=3.4 in]{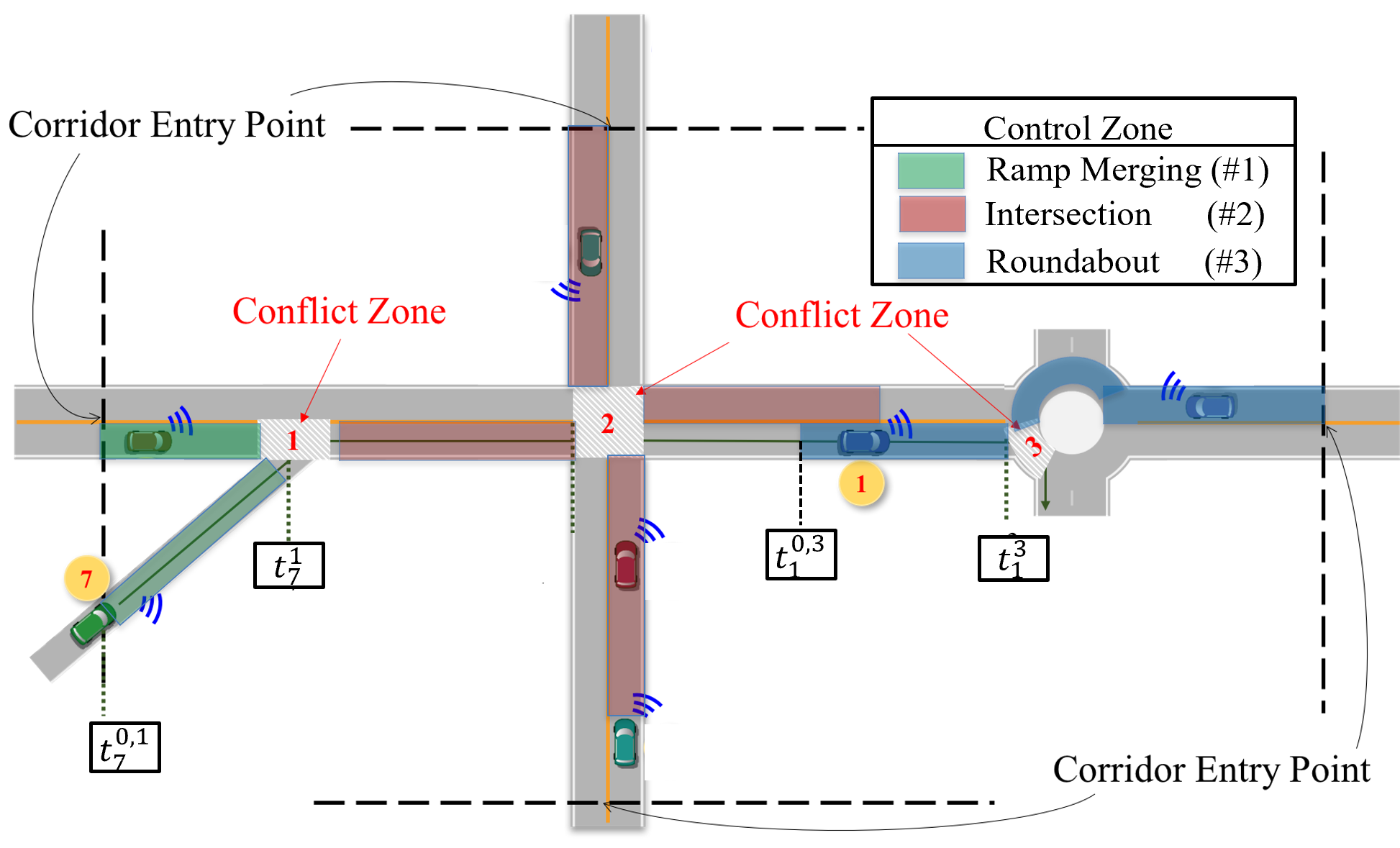} \caption{Corridor with
connected and automated vehicles.}%
\label{fig:corridor}%
\end{figure}

We consider a corridor (Fig. \ref{fig:corridor}) that consists of several conflict zones (e.g., a merging area, an intersection, and a roundabout), where potential lateral collision of vehicles may occur. Upstream of each conflict zone, we define a \textit{control zone}, inside of which, the vehicles can communicate with each other. The dimension of each control zone is restricted by the communication range of an associated \textit{coordinator}, which records the vehicle queue inside the control zone. Note that the coordinator is not involved in any decision on the CAV operation and only serves to coordinate information with the CAVs. The communication range of the coordinator can be adjustable and its length could be extended as needed. For clarity, we illustrate the boundary of the corridor as indicated by dashed lines and the limits of each control zone by shaded rectangles (Fig.~\ref{fig:corridor}). Note that we only coordinate CAVs inside the control zone of each conflict zone.

We consider a corridor with a set of conflict zones $\mathcal{Z} \subset \mathbb{N}$. Let $\mathcal{N}_z(t)=\{1,2,\dots,N_z(t)\}$ be a queue of CAVs inside the control zone of a conflict zone $z$, where $N_z(t)\in\mathbb{N}$ is the number of CAVs in the control zone of $z$ at time $t\in\mathbb{R}^{+}$. When a CAV enters the control zone, it broadcasts its route information to the coordinator of this conflict zone. The coordinator then assigns a unique integer $i\in\mathbb{N}$ that serves as the identification of CAVs inside the corridor. Let $t_i^{0,z}$ be the initial time that CAV $i$ enters the control zone of $z \in \mathcal{Z}$, and $t_{i}^{z}$ be the time for CAV $i$ that enters $z$. For example, for CAV $\# 7$ (Fig.~\ref{fig:corridor}), $t_{7}^{0,1}$ is the time that it enters the control zone of conflict zone $\#1$, which is also the time that it enters the corridor, and $t_{7}^{1}$ is the time that it enters the conflict zone $\#1$. Similarly, CAV $\# 1$ enters the control zone of conflict zone $\#3$ at $t_{1}^{0,3}$, and enters the conflict zone $\#3$ at $t_{1}^{3}$. The CAV index given by the coordinator is removed from the queue $\mathcal{N}_z$ once the vehicle $i$ exits the conflict zone $z$.

To avoid any possible lateral collision, there is a number of ways to compute $t_{i}^{z}$ for each CAV $i$. In what follows, we present a decentralized framework in which we formulate an upper-level optimal control problem for determining the time $t_{i}^{z}$ that each CAV $i$ will enter the conflict zone $z \in \mathcal{Z}$, and then address a lower-level control problem that will yield for each CAV the optimal control input (acceleration/deceleration) to achieve the assigned time $t_{i}^{z}$ (upon arrival of CAV $i$) without collision.

\subsection{Vehicle model, Constraints, and Assumptions}
Each CAV $i\in\mathcal{N}_z(t)$ is modeled by a second order dynamics
\begin{equation}%
\dot{p}_{i} = v_{i}(t), \quad \dot{v}_{i} = u_{i}(t),
\label{eq:model2}
\end{equation}
where $p_{i}(t)\in\mathcal{P}_{i}$, $v_{i}(t)\in\mathcal{V}_{i}$, and $u_{i}(t)\in\mathcal{U}_{i}$ denote the position, speed and control input (acceleration/deceleration) of each CAV $i$ in the control zone. Let $x_{i}(t)=\left[p_{i}(t) ~ v_{i}(t)\right]  ^{T}$ denote the state of each CAV $i$, with initial value at the entry of the control zone of conflict zone $z\in \mathcal{Z}$ given as $x_{i}^{0,z}=\left[p_{i}^{0,z}  ~ v_{i}^{0,z}\right]  ^{T}$, where $p_{i}^{0,z}= p_{i}(t_{i}^{0,z})$ and $v_{i}^{0,z}= v_{i}(t_{i}^{0,z})$. 
To ensure that each vehicle control input and speed are within a given admissible range,  we impose the following constraints
\begin{equation}%
\begin{split}
u_{i,min} &  \leq u_{i}(t)\leq u_{i,max},\quad\text{and}\\
0 &  \leq v_{min}\leq v_{i}(t)\leq v_{max},\quad\forall t\in\lbrack t_{i}%
^{0,z},t_{i}^{z}],
\end{split}
\label{speed_accel_constraints}%
\end{equation}
where $u_{i,min}$, $u_{i,max}$ are the minimum deceleration and maximum
acceleration for each CAV $i\in\mathcal{N}_z(t)$, and $v_{min}$, $v_{max}$ are the minimum and maximum speed limits respectively. For simplicity, we do not consider vehicle diversity, and thus we set $u_{i,min}=u_{min}$ and $u_{i,max}=u_{max}$.
For each CAV $i\in \mathcal{N}_z(t)$, the lateral collision is possible within the set $\Gamma_i \overset{\Delta}{=} \{t\,\,|\,t\in  [t_i^{z}, t_i^{z}+\rho]\}$,
where, $\rho$ is the safety time headway to avoid lateral collision. Lateral collision between any two CAVs $i,j\in \mathcal{N}_z(t)$ can be avoided if the following constraint holds,
\begin{equation}
\Gamma_i \cap \Gamma_j=\varnothing, ~  \forall t\in [t_i^{z}, t_i^{z}+\rho], \quad i,j\in \mathcal{N}_z(t).
\label{eq:lateral_constraint}
\end{equation} 
To ensure the absence of rear-end collision of two consecutive CAVs traveling on the same lane, we impose the following condition
\begin{equation}
\begin{split}
s_{i}(t)= p_{k}(t)-p_{i}(t) \ge \delta(t),~ \forall t\in [t_i^{0,z}, t_i^z].
\label{eq:rearend}
\end{split}
\end{equation}
Here, $s_{i}(t)$ denotes the distance between CAV $i$ and CAV $k$ which is physically immediately ahead of $i$. The minimum safe distance $\delta(t)$ is a function of speed $v_i(t)$. Since we consider an urban traffic corridor, the average speed does not exhibit significant variations. Therefore, we can
consider that the safe distance $\delta(t)=\delta$ is constant.
In the  modeling framework described above, we impose the following assumptions:
\begin{assum} \label{ass1}
All vehicles are connected and automated, i.e., 100\% penetration rate of CAVs.
\end{assum}

\begin{assum} \label{ass2}
For each CAV, none of the constraints \eqref{speed_accel_constraints} and \eqref{eq:rearend} is active at $t^{0,z}_i$.
\end{assum}

\begin{assum}\label{ass3}
Each CAV $i$ has proximity sensors and can measure local information without errors or delays.
\end{assum}

\begin{assum} \label{ass: no_lane_change}
The corridor only contains single-lane road segments. The vehicles traveling in the corridor do not change lanes except to make necessary turns.
\end{assum}

\begin{assum}\label{ass:constant_MZ_speed}
The speed of each CAV $i$ inside the conflict zone is constant.
\end{assum}

The first assumption limits the scope of our paper to the control of CAVs in an idealized environment where all vehicles are automated and connected to each other. Addressing different penetration of CAVs is the object of ongoing work. The second assumption ensures that the solution of the optimal control problem starts from a feasible state and control input. 
The third assumption might impose barriers in a potential deployment of the proposed framework. However, we could extend our results in the case that this assumption is relaxed, if the noise in the measurements and delays are bounded. 
The fourth assumption simplifies the upper-level optimal control problem so as to avoid implications related to lane changing. 
Finally, the last assumption is imposed to enhance safety awareness. However, it could be modified appropriately, if necessary, as discussed by \cite{Malikopoulos2017}.

\begin{defn} \label{def: subsets}
Let $i-1, i \in \mathcal{N}_z(t)$ be two CAVs inside the control zone traveling towards the corresponding conflict zone $z$. Depending on the physical location and trajectory inside the control zone with respect to CAV $i$, CAV $i-1$ belongs to one of the following three subsets of $\mathcal{N}_z(t)$ with respect to CAV $i \in \mathcal{N}_z(t)$: 
\begin{enumerate}
\item{$\mathcal{R}^{z}_i$ contains all CAVs that travel in the same lane with CAV $i$ towards the conflict zone $z$, having travel paths that can cause rear-end collision.}

\item{$\mathcal{C}^{z}_i$ contains all CAVs from different roads   having travel paths that can cause lateral collision with CAV $i$ in conflict zone $z$.}

\item{$\mathcal{O}^{z}_i$ contains all CAVs from different roads having travel paths that cannot cause lateral or rear-end collision with CAV $i$ in conflict zone $z$}
\end{enumerate}
\end{defn}
Upon arrival at the entry of the control zone of conflict zone $z\in \mathcal{Z}$ at time $t_i^{0,z}$, CAV $i\in\mathcal{N}_z(t)$ needs to compute the time $t_i^{z}$. In general, a value of $t_i^{z}$ that satisfies the safety constraints \eqref{eq:lateral_constraint} and \eqref{eq:rearend} may depend on the preceding CAV in the control zone. Next, we address the question of identifying the appropriate $t_i^{z}$ for each CAV through an upper-level optimization problem.

\subsection{Upper-level optimization problem}
To fully utilize the capacity of the traffic network of CAVs, we formulate a throughput maximization problem, in terms of minimizing the gaps between CAVs in the conflict zones, i.e., minimizing the total time to process CAVs in the network, subject to the constraints \eqref{speed_accel_constraints}, \eqref{eq:lateral_constraint} and \eqref{eq:rearend}. Note that for $i=1$, the safety constraint is not active since there is no prior CAV in the control zone, which implies that $t_1^{z}$ is not constrained and can be determined outside the optimization framework. Thus, for each control zone of conflict zone $z\in \mathcal{Z}$, we formulate the following optimization problem:
\begin{gather} \label{eq:throughput}
\min_{\textbf{t}_{(2:N_z(t))}} \sum_{i=2}^{{N}_z(t)}(t_i^{z} - t_{i-1}^{z}) = \min_{\textbf{t}_{(N_z(t))}} (t_{{N}_z(t)}^{z} - t_1^{z}), \\
\text{subject to}: \eqref{eq:model2}, \eqref{speed_accel_constraints}, \eqref{eq:lateral_constraint}, \eqref{eq:rearend}, \nonumber 
\end{gather}
where, $_{\textbf{t}_{(2:N_z(t))}}= [t_2^z, \dots , {t}_{N_z(t)}^z]$.
The solution of \eqref{eq:throughput} yields the optimal time ${t_i^{z}}^*, i\geq 2, z\in \mathcal{Z}$, which designates the entry times of each CAV in each conflict zone so as to maximize the throughput of the corresponding bottleneck.

In what follows, we discuss how the lateral collision safety constraint is addressed in the solution of \eqref{eq:throughput}. We also show that the solution has an iterative structure and depends only on the state and control constraint \eqref{speed_accel_constraints} as well as the safety constraint \eqref{eq:rearend}. 
To obtain the optimal solution of ${t_i^{z}}^*$ for CAV $i\in\mathcal{N}_z(t)$ at the conflict zone $z\in\mathcal{Z}$, we first consider the case when $\mathcal{C}_i^z(t)$ is empty, thus the entry time for CAV $i$ at $z$ depends only on some CAV $k = (i-1) \in \mathcal{R}^{z}_i$, where $k$ is physically immediately ahead of $i$ on the roadway segment inside the control zone. In this case, the minimum time ${t_i^{z}}$ for CAV $i$ to enter the conflict zone $z$ is designated by the rear-end safety constraint \eqref{eq:rearend}, and in particular, by the safe headway, $\rho,$ that a CAV $i$ should maintain while following CAV $k$, i.e., $t_i^{z}=t_{k}^{z}+\rho.$
In this context, we need to find the bound of $t_i^{z}$ to ensure feasibility of the solution. Consider the maximum and minimum speeds that CAV $i$ could achieve. The value of $t_i^{z}$ is then given by
\begin{align} \label{eq: t_i^m_j_2}
{t_i^{z}}=\max\big\{\min\{t_{k}^{z}+\rho,~ t^{z,max}_{i}\}, t^{z,min}_{i} \big\}, 
\end{align}
where $t^{z,max}_{i}$ and $t^{z,min}_{i}$ is the longest and shortest possible travel time of a CAV $i$ between the entry and exit of the control zone of the  conflict zone $z \in \mathcal{Z}$ corresponding to the minimum, $v_{min},$ and maximum, $v_{max},$ speed limit respectively.
Note that condition \eqref{eq: t_i^m_j_2} ensures that the time $t_i^{z}$ that CAV $i$ will be entering the conflict zone $z$ is feasible and can be attained based on the imposed speed limits in the corridor. From \eqref{eq: t_i^m_j_2}, the safety constraint between CAVs traveling in the same lane is guaranteed at ${t_i^{z}}$. 
We now turn our attention to the case where possible lateral collision might occur if $\mathcal{C}^{z}_i(t)$ is non-empty. In this case, the minimum time $t_i^{z*}$ for CAV $i$ to enter the conflict zone $z$ is constrained by both the lateral collision \eqref{eq:lateral_constraint} and rear-end collision \eqref{eq:rearend} constraints.
\begin{defn} \label{def_a}
We define the set $\mathcal{A}^{z}_i\subset\mathcal{C}^{z}_i$, $\mathcal{A}^{z}_i \mathrel{\mathop:}= \{j \in \mathcal{C}^{z}_i~|~t_j^{z} \geq t_i^{z}\},$ that includes any CAV $j \in \mathcal{C}^{z}_i$ whose entry time at conflict zone $z$ is later than $t_i^{z}$, and the set $\mathcal{L}^{z}_i\subset\mathcal{A}^{z}_i$, $\mathcal{L}^{z}_i \mathrel{\mathop:}= \{j \in \mathcal{A}^{z}_i~|~t_j^{z}+\rho \leq t_{j+1}^{z}-\rho \}$, that includes any CAV $j \in \mathcal{C}^{z}_i$ whose entry time satisfy \eqref{eq:rearend}.
\end{defn}
Considering possible lateral collisions at conflict zone $z$, for all $z \in \mathcal{Z}$, we obtain the following result.
\begin{thm} \label{theo: time}
The solution ${t_i^{z}}^*, \forall i \geq 2, \forall z \in \mathcal{Z},$ of  \eqref{eq:throughput} is recursively determined through
\begin{equation} \label{eq: t_i^m_j_*}
\resizebox{0.5\textwidth}{!}{%
${t_i^{z}}^* = 
\begin{cases}
\max\big\{\max\{t_c^{z}\} + \rho_{i}, t_i^{z} \big\}, & \text{if } \forall c \in \mathcal{C}^{z}_i \text{  and } \nexists ~c \in \mathcal{A}^{z}_i,\\
t_i^{z}, &\text{if } ~\exists ~a \in \mathcal{A}^{z}_i \text{ and } t_i^{z}+\rho \leq \min\{t_a^{z}\},\\
\min\{t_b^{z}\} + \rho_{i}, & \text{if }~ \exists ~b \in \mathcal{L}^{z}_i \text{ and } t_i^{z}+\rho_{i} > \min\{t_a^{z}\}, \\ 
\max\{t_a^{z}\} + \rho_{i}, &\text{if }~ \exists ~a \in \mathcal{A}^{z}_i \text{ and } \nexists ~a \in \mathcal{L}^{z}_i.
\end{cases}$
}
\end{equation}
\end{thm}
\textbf{PROOF.}
Based on Definition \ref{def_a}, there are three cases to consider for ${t_i^{z}}^*$.\newline
\textit{Case 1}: If $\mathcal{A}^{z}_i\neq\emptyset$, all CAVs in $\mathcal{C}^{z}_i$ will be entering the conflict zone $z$ earlier than $t_i^{z}$, which, by Definition \ref{def_a}, implies ${t_i^{z}}^* =\bigg\{ \begin{matrix}
t_i^{z}, \\
\max\{t_c^{z}\} + \rho,
\end{matrix} \quad \quad \forall c \in \mathcal{C}^{z}_i.$

Hence, we have ${t_i^{z}}^* = \max\big\{\max\{t_c^{z}\} + \rho, t_i^{z} \big\}, \forall c \in \mathcal{C}^{z}_i$.
\newline
\textit{Case 2}: If $\mathcal{A}^{z}_i\neq\emptyset$ and $\mathcal{L}^{z}_i\neq\emptyset$, we consider two cases: (i) if the earliest entry time of CAVs in the set $\mathcal{A}^{z}_i$ is later than $t_i^{z}$ plus a safe headway $\rho$, then the minimum entry time of CAV $i$ is $t_i^{z}$, which satisfies the safety constraints to avoid both lateral and rear-end collisions; (ii) the optimal value of $t_i^{z}$ is the earliest possible time slot between the entry times of two consecutive CAVs in the set $\mathcal{L}^{z}_i$. Hence, we have
\begin{align} \label{eq: t_i^m_j_x_2}
{t_i^{z}}^* = 
\begin{cases}
t_i^{z}, & \text{if } t_i^{z}+\rho \leq \min\{t_a^{z}\},\\
\min\{t_b^{z}\} + \rho, & \text{if }  t_i^{z}+\rho > \min\{t_a^{z}\},
\end{cases} \nonumber \\
\forall a \in \mathcal{A}^{z}_i, \forall b \in \mathcal{L}^{z}_i.
\end{align}
\textit{Case 3}: Finally, if $\mathcal{L}^{z}_i\neq\emptyset$, this implies that there is no available time slot between the entry times of two CAVs in $\mathcal{A}^{z}_i$. In this case, CAV $i$ will be entering the conflict zone after the last CAV in $\mathcal{A}^{z}_i$ to avoid lateral collision, which implies ${t_i^{z}}^* = \max\{t_a^{z}\} + \rho, \forall a \in \mathcal{A}^{z}_i,$ if there not exist $~a \in \mathcal{L}^{z}_i$.
Combining the above results, we obtain ${t_i^{z}}^*$ in \eqref{eq: t_i^m_j_*}, which completes the proof.
 
Theorem \ref{theo: time} yields the sequence that the CAVs will be traveling through each control zone. Each CAV $i$ follows the above policy to determine the time ${t_i^{z}}^*$ that it will be entering the conflict zone $z \in \mathcal{Z}$ upon arrival at the entry of the control zone. Once the entry time ${t_i^{z}}^*$ is computed, it is stored in the coordinator and it is not changed. Thus, the next CAV $i+1$, upon its arrival at the entry of the control zone, will search for feasible times to cross the conflict zone based on the available time slots. The recursion is initialized when the first CAV enters the boundary of the corridor, i.e., it is assigned $i = 1$. In this case, $t_1^{z}, \forall z\in \mathcal{Z},$ can be externally assigned as the desired entry time of this CAV whose behavior is unconstrained.

\subsection{Low-level optimal control problem} \label{sec:sol}
When a CAV $i$ enters the corridor, it communicates with the other CAVs (Assumption \ref{ass1}) and the coordinator broadcasts information without any errors or delays (Assumption \ref{ass3}). The coordinator assigns a unique identity to each CAV and receives back some information at the time each CAV arrives at the entry of the corridor, as defined next. 

\begin{defn}
For each CAV $i$, we define the \textit{information set} ${Y}_i(t)$ as 
\begin{equation}
{Y}_i(t) \triangleq \{ p_i(t), v_i(t), o_i, {t_i^{z}}^*\}, z\in\mathcal{Z}, t\in [t_i^0, t_i^z],
\end{equation}
where $p_i(t)$, $v_i(t)$ are the position and speed of CAV $i$ inside
the corridor, $o_i$ is the route that CAV $i$ travels inside the corridor, and ${t_i^{z}}^*$ is the time for CAV $i$ to enter the conflict zone $z$ given by \eqref{eq: t_i^m_j_*}. 
\end{defn}

As discussed in the previous section, the time ${t_i^{z}}^*$ for CAV $i$ is determined in a recursive manner based on the information received from the coordinator. Therefore, once CAV $i$ enters each of the control zones, immediately all information in ${Y}_i(t)$ becomes available for $i$ and is stored in the coordinator accessible for next arriving CAV $i+1$.

In the low-level optimal control problem, the objective is to minimize the control input (acceleration/deceleration) for each CAV $i\in\mathcal{N}_z(t)$ from the time $t_i^{0,z}$
that $i$ enters the control zone until the time $t_i^{z}$ that it exits the control zone under the hard safety constraint to avoid rear-end collision. By minimizing each CAV's acceleration/deceleration, we minimize transient engine operation. Thus, we can have direct benefits in fuel consumption and emissions since internal combustion engines are optimized over steady state operating points (constant torque and speed).
Therefore, the optimization problem for each CAV $i\in\mathcal{N}_z(t)$ is to minimize the $L^2$-norm of the control input in $[t_i^{0,z}, t_i^z]$, formulated as follows:
\begin{gather}\label{eq:decentral}
\min_{u_i(t)\in U_i}J_{i}(u(t))= \min_{u_i(t)\in U_i} \frac{1}{2} \int_{t^{0,z}_i}^{t^z_i} u^2_i(t)~dt,\\ 
\text{subject to}%
:\eqref{eq:model2}, \eqref{speed_accel_constraints},\eqref{eq:rearend},\nonumber\\
\text{given }t_{i}^{0,z},~ p_{i}(t_{i}^{0,z})=0, ~ v_{i}^{0,z},~{t_i^{z}}^*, \text{ and } p_{i}({t_i^{z}}^*)=p_z.\nonumber
\end{gather}
Note that we do not include the lateral collision constraint \eqref{eq:lateral_constraint} in \eqref{eq:decentral}, since it has been addressed in the upper-level optimization problem. On the contrary, we explicitly include the rear-end safety constraint. The  problem formulation with the state and control constraints requires the constrained and unconstrained arcs of the state and control input to be pieced together to satisfy the Euler-Lagrange equations and necessary condition of optimality. Let $\textbf{S}(t,\textbf{x}_i(t),u_i(t))=[v_{i}(t) - v_{max}~~~
v_{min}(t) - v_{i}(t)~~~
p_{i}(t) - p_{k}(t)-\delta
]^T$ be the vector of constraints that are not explicit functions of the control input $u_i(t)$.
We take successive total time derivatives of $\textbf{S}(t,\textbf{x}_i(t),u_i(t))$ until we obtain an expression that is explicitly dependent on $u_i(t)$. If $q$ time derivatives are required for a specific constraint of $\textbf{S}_i(t,\textbf{x}_i(t),u_i(t))$, we refer to that constraint as a $q$th-order state variable inequality constraint; see \cite{bryson1975applied}. Note that we have 1st-order speed constraint and 2nd-order rear-end safety constraint in $\textbf{S}_i(t,\textbf{x}_i(t),u_i(t))$. The 2nd-order rear-end safety constraint plays the role of a control variable constraint on the constrained arc,
\begin{equation}
    S_i^{(2)}(x_i(t),u_i(t),t) = u_i(t) - u_k(t) = 0. \label{eq:uk_constraint}
\end{equation}
From \eqref{eq:decentral}, the CAV dynamics \eqref{eq:model2}, state and control constraints \eqref{speed_accel_constraints}, and the rear-end safety constraint \eqref{eq:rearend} for each CAV $i\in\mathcal{N}_z(t)$, we formulate the Hamiltonian function 
\begin{gather}
H_{i}\big(t, p_{i}(t), v_{i}(t), u_{i}(t)\big)
=\frac{1}{2} u(t)^{2}_{i} + \lambda^{p}_{i} \cdot v_{i}(t)  + \lambda^{v}_{i} \cdot u_{i}(t) \nonumber\\
+ \eta^{a}_{i} \cdot u_i(t) - \eta^{b}_{i} \cdot u_i(t) + \eta^{c}_{i} \cdot (u_i(t)-u_k(t)) \nonumber \\ +  \mu^{d}_{i} \cdot (u_i(t)-u_{max})+\mu^{e}_{i} \cdot (u_{min}-u_i(t)) ,
\label{eq:lagrange-hamil}
\end{gather}
where $\lambda^{p}_{i},~\lambda^{v}_{i}$ are the co-state components, and $\eta^{a}_{i}, \eta^{b}_{i}, \eta^{c}_{i}, \mu^{d}_i \mu^{e}_i$ are the Lagrange multipliers satisfying the complimentary slackness conditions based on the state and control constraints in \eqref{speed_accel_constraints} and \eqref{eq:rearend}.
The Euler-Lagrange equations become
\begin{gather}\label{eq:ELp}
\dot\lambda^{p}_{i}(t) = - \frac{\partial H_i}{\partial p_{i}} = 0, \quad \dot\lambda^{v}_{i}(t) = - \frac{\partial H_i}{\partial v_{i}} = -\lambda^{p}_{i}.
\end{gather}
The necessary condition for optimality is
\begin{equation}
\label{eq:KKT1}\frac{\partial H_i}{\partial u_{i}} = u_{i}(t) + \lambda^{v}_{i}+ \eta_i^a - \eta_i^b + \eta_i^c +\mu_i^d - \mu_i^e = 0.
\end{equation}
\subsubsection{State and control constraint is not active} \label{case1}
When the inequality state and control constraints are not active, $ \eta_i^a = \eta_i^b = \eta_i^c =\mu_i^d = \mu_i^e = 0$, applying the necessary condition \eqref{eq:KKT1}, the optimal control is
\begin{equation}
u_i(t) + \lambda^{v}_{i}= 0, \quad i \in\mathcal{N}_z(t). \label{eq:23}
\end{equation}
From (\ref{eq:ELp}) we have $\lambda^{p}_{i}(t) = a_{i}$, and $\lambda^{v}_{i}(t) = -(a_{i}\cdot t + b_i)$. 
The coefficients $a_{i}$, $b_{i}$
are constants of integration corresponding to each CAV $i$. From \eqref{eq:23} the optimal control input (acceleration/deceleration), speed and position as a function of time are given by
\begin{gather}
u^{*}_{i}(t) = (a_{i} \cdot t + b_i), ~ \forall t \ge t^{0,z}_{i}, \label{eq:24}\\
v^{*}_{i}(t) = \frac{1}{2} a_{i} \cdot t^2 + b_{i} \cdot t +c_{i}, ~ \forall t \ge t^{0,z}_{i},\label{eq:25}\\
p^{*}_{i}(t) = \frac{1}{6} a_{i} \cdot t^3 +\frac{1}{2} b_{i} \cdot t^2 + c_{i}\cdot t +d_{i}, ~ \forall t \ge t^{0,z}_{i}, \label{eq:26}%
\end{gather}
where $c_{i}$ and $d_{i}$ are constants of integration which can be computed at each time $t, t^{0}_{i} \le t \le t^{z}_{i}$, using the values of the control input, speed, and position of each CAV $i$ at $t$, the position $p_{i}(t^{z}_{i})$, and $\lambda^{v}_{i}(t^{z}_{i}) =0$. \newline
To derive the constrained solution of \eqref{eq:decentral}, we first start with the unconstrained arc and derive the solution using \eqref{eq:24}. If the solution violates any of the constraints, then we re-solve the problem with the constrained and unconstrained arcs pieced together. 
This process is repeated until the solution does not violate any other constraints. 
The simple nature of the optimal control and states
in \eqref{eq:24} through \eqref{eq:26} makes the online solution of \eqref{eq:decentral} computationally
feasible, even with the additional burden of checking
for active constraints.
In what follows, we address the optimization problem \eqref{eq:decentral} with the activation of the constrained case corresponding to the rear-end collision only. The other constrained cases related to the state, i.e., speed, $v_i$, and control, $u_i$ as in \eqref{speed_accel_constraints}, are similar to the cases presented by \cite{Malikopoulos2017} and  \cite{Mahbub2020ACC-1}, and thus are omitted.

\subsubsection{Rear-end safety constraint is active}\label{case2}
Suppose a CAV starts from a feasible state and control at $t=t_{i}^0$ and at some time $t=t_{1}$, the rear-end safety constraint $(p_i(t)-p_k(t)+\delta) \le 0$ is activated. In this case, $\eta_i^a = \eta_i^b =\mu_i^d = \mu_i^e = 0$.
From \eqref{eq:KKT1}, the
optimal control is given by 
\begin{gather}\label{eq:case3u}
u_{i}(t) + \lambda^{v}_{i} + \eta_i^c(t)= 0, ~\forall t \ge t_{1}.
\end{gather}
The state trajectory entering onto the 2nd-order rear-end safety constraint boundary must satisfy the following tangency conditions
\begin{align}
N(x_i(t),t) = \begin{bmatrix}
p_i(t)-p_k(t)+\delta\\ 
v_i(t) - v_k(t)
\end{bmatrix} = 0,\label{eq:tangency}
\end{align}
where, $N(x_i(t),t)$ is the $q$-component vector function of the 2nd order safety tangency constraints. The tangency constraints in \eqref{eq:tangency} also apply to the state trajectory at the exit of the constraint arc. Since the optimal solution of the preceding CAV $k\in \mathcal{N}_z(t)$ is known a priori, from \eqref{eq:uk_constraint} and \eqref{eq:tangency}, the optimal solution for CAV $i\in \mathcal{N}_z(t)$ in the constrained arc is derived from
$\textbf{S}_i(t,\textbf{x}_i(t),u_i(t))=0$ and is
\begin{gather}
    u_i^*(t) = u_k^*(t),~ v_i^*(t) = v_k^*(t),~ \text{and}~ p_i^*(t) = p_k^*(t)-\delta.\label{eq:optimal} 
\end{gather}
The equations in \eqref{eq:tangency} form a set of interior boundary conditions where the co-states $\lambda_i^p(t)$ and $\lambda_i^v(t)$ in general have discontinuity at the junction points, i.e., entry and exit points of the state trajectory between the constrained and the unconstrained arcs. However, the control trajectory may or may not have discontinuities at the junction points. Next, we address the jump conditions at the entry junction point $t=t_1$.
At time $t = t_{1}$, when the safety constraint is activated, we have a junction point between the unconstrained and constrained arcs. Let $t^{-}_{1}$ represents the time instance just before $t_{1}$, and $t^{+}_{1}$ signifies just after $t_{1}$. The state trajectories are continuous at the junction points. Thus, we have
\begin{gather}
    p_i(t^{-}_{1}) = p_i(t^{+}_{1}),  \quad   v_i(t^{-}_{1}) = v_i(t^{+}_{1}). \label{eq:speed-continuity}
\end{gather}
The jump conditions at $t_{1}$ can be written as
\begin{gather}
{\lambda_i^p}(t^{-}_{1}) = {\lambda_i^p}(t^{+}_{1}) + \pi^T\cdot\frac{\partial N(x_i(t),t)}{\partial p_i(t)}\bigg{|}_{t=t1}, \label{eq:jump1}\\
{\lambda_i^v}(t^{-}_{1}) = {\lambda_i^v}(t^{+}_{1}) + \pi^T\cdot\frac{\partial  N(x_i(t),t)}{\partial v_i(t)}\bigg{|}_{t=t1}, \label{eq:jump2}\\
H(t^{-}_{1}) = H(t^{+}_{1}) -  \pi^T\cdot\frac{\partial  N(x_i(t),t)}{\partial t}\bigg{|}_{t=t1}, \label{eq:jump3}\\
\frac{\partial H(t^{-}_{1})}{u_i(t)}  =\frac{\partial H(t^{+}_{1})}{u_i(t)}. \label{eq:jump4}
\end{gather}
In \eqref{eq:jump1}-\eqref{eq:jump4}, $\pi^T = [\pi_1 ~ \pi_2]$ is a vector of constant Langrange multipliers to be determined so that the condition in \eqref{eq:tangency} is satisfied.
From \eqref{eq:jump3}, using \eqref{eq:optimal}-\eqref{eq:jump4}, we obtain $\frac{1}{2} (u_{i}(t^{-}_{1})^{2}-u_{i}(t^{+}_{1})^{2}) + \lambda_i^v(t^{+}_{1}) \cdot (u_{i}(t^{-}_{1})-u_{i}(t^{+}_{1}))=0$. This yields two cases: either $ (u_{i}(t^{-}_{1})-u_{i}(t^{+}_{1}))=0$ or $\frac{1}{2} (u_{i}(t^{-}_{1})+u_{i}(t^{+}_{1}))+\lambda_i^v(t^{+}_{1})=0$. Both conditions lead to $u_{i}(t^{-}_{1})=u_{i}(t^{+}_{1})$,
which indicates that the control trajectory is continuous at the entry junction point at $t=t_1$.
Finally, using the continuity in control and \eqref{eq:jump2} in \eqref{eq:jump4}, we obtain $\eta_i^c(t_1^+) = \pi_2$. With two junction points at time $t=t_1$ and $t=t_2$, we have a constrained arc between two unconstrained arcs. Since we have multiple arcs pieced together at the junction points, we differentiate the constants of integration for the state and control trajectory by adding a superscript $h$ representing the order of appearance of the arcs. Therefore, we represent the constants of integration as  $(a_i^h,b_i^h,c_i^h,d_i^h)$, where $h=1,2$ corresponds to the first and the last unconstrained arc respectively. The control trajectory of CAV $i$ considering the constrained and unconstrained arcs can be written as,
\begin{equation}\label{eq:control-solution}
{u}_i^*(t) = \left\{
\begin{array}
[c]{ll}%
& \mbox{$a_i^1\cdot t + b_i^1, ~~ t_i^{0,z}\ge t \le t_1$},\\
& \mbox{$u_k^*(t), ~~~~~~~~ t_1< t < t_2$},\\
& \mbox{$a_i^3\cdot t + b_i^3, ~~ t_2\ge t \le t_i^z$}.\\
\end{array}
\right.
\end{equation}
The constants of integration, along with the junction points $t_1$ and $t_2$ can be computed by solving \eqref{eq:24}-\eqref{eq:26} and \eqref{eq:optimal}-\eqref{eq:control-solution} with appropriate initial, boundary and transversality conditions.
\section{Simulation Results} \label{sec:sim}

\begin{figure}
    \centering
    \includegraphics[scale=0.52]{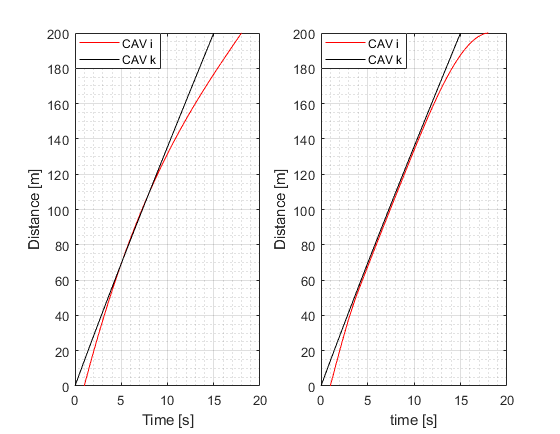}
    \caption{Unconstrained (left) and rear-end safety constrained state trajectory (right) of CAV i with respect to its immediately preceding vehicle k.}
    \label{fig:safety}
\end{figure}

To validate the effectiveness of the rear-end safety constrained formulation, we present two cases in Fig. \ref{fig:safety}, where a leading CAV $k$ and a following CAV $i$ are both cruising with the optimal control input. In the left panel of Fig. \ref{fig:safety}, the following CAV $i$ derives its control input according to \eqref{eq:24}, and activates the rear-end safety constraint with respect to the immediately preceding CAV $k$ with two junction points. In the right panel of Fig. \ref{fig:safety}, CAV $i$ derives its control input using \eqref{eq:control-solution} subject to safety constrained optimization.

To validate the proposed approach for multiple traffic scenarios, we use a simulation network of Mcity created in PTV VISSIM environment. We define a corridor consisting of four conflict zones: (1) a merging roadway, (2) a speed reduction zone, (3) a roundabout, and (4) an intersection. 
To evaluate the network performance with the proposed control framework, we define two scenarios as follows:\newline
\textbf{Scenario 1:} baseline, i.e., 0\% CAV penetration rate. All vehicles in the network are non-connected and non-automated. In this case, the Wiedemann car following model \cite{wiedemann1974} built in VISSIM is applied. 1.2 $s$ time headway is adopted to estimate the minimum allowable following distance.\newline
\textbf{Scenario 2:} optimal control, i.e., 100\% CAV penetration rate. The proposed control framework is integrated to generate the optimal acceleration/deceleration profile for each CAV in the network. 

\begin{figure} [!ht]
	\centering
	\includegraphics[width=0.42\textwidth]{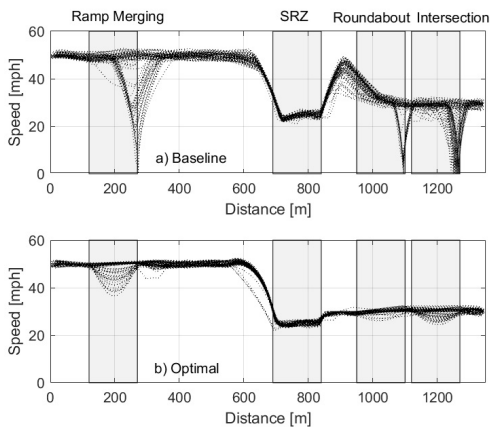}   
	\caption{Vehicle trajectories inside the corridor for (a) baseline and (b) optimal controlled case. The control zone for each of the conflict zones are shown for comparison.} 
	\label{fig:trajectory}
\end{figure}

The CAV speed trajectories under 0\% and 100\% CAV penetration rate in the corridor are illustrated in Fig. \ref{fig:trajectory}. In the baseline scenario with 0\% CAV penetration rate, CAVs traveling along the corridor need to yield to mainline traffic, and wait in the signalized intersection. Thus, we observe high fluctuations in their speed profiles under the baseline scenario at the proximity of the conflict zones (see upper panel of Fig. \ref{fig:trajectory}). In the optimal control scenario under 100\% CAV penetration rate, CAVs travel through the corridor without stop-and-go driving (see lower panel of Fig. \ref{fig:trajectory}). The latter enables  CAVs to have smoother speed trajectory affecting the uncontrolled upstream and downstream area of the control zone.  We observe $9\%$ improvement in terms of travel time and an average of 47\% savings in total fuel consumption in the optimal control scenario compared to the baseline one.  

\section{Concluding Remarks and Discussion} \label{sec:conc}
In this paper, we investigated the optimal coordination of CAVs in a corridor. We presented a two-level optimization problem in which we maximize traffic throughput in the upper-level problem, and derive a closed-form analytical solution that yields the optimal control input for each CAV, in terms of fuel consumption, in the low-level problem.
We derived a closed-form analytical solution that considers safety constraints. We showed through simulation that vehicle coordination can reduce stop-and-go driving improve efficiency in the corridor under hard safety constraints. 
One limitation of the proposed approach is that it does not consider the interaction of the conflict zones which might have significant implications as the traffic volume increases. Future research should address address this interdependence as it could enhance the benefits in energy and travel time significantly.


\bibliographystyle{plainnat_automa}
\bibliography{corridor_ref}
\par
\begin{wrapfigure}{L}{0.17\textwidth}
\centering
\includegraphics[width=0.21\textwidth]{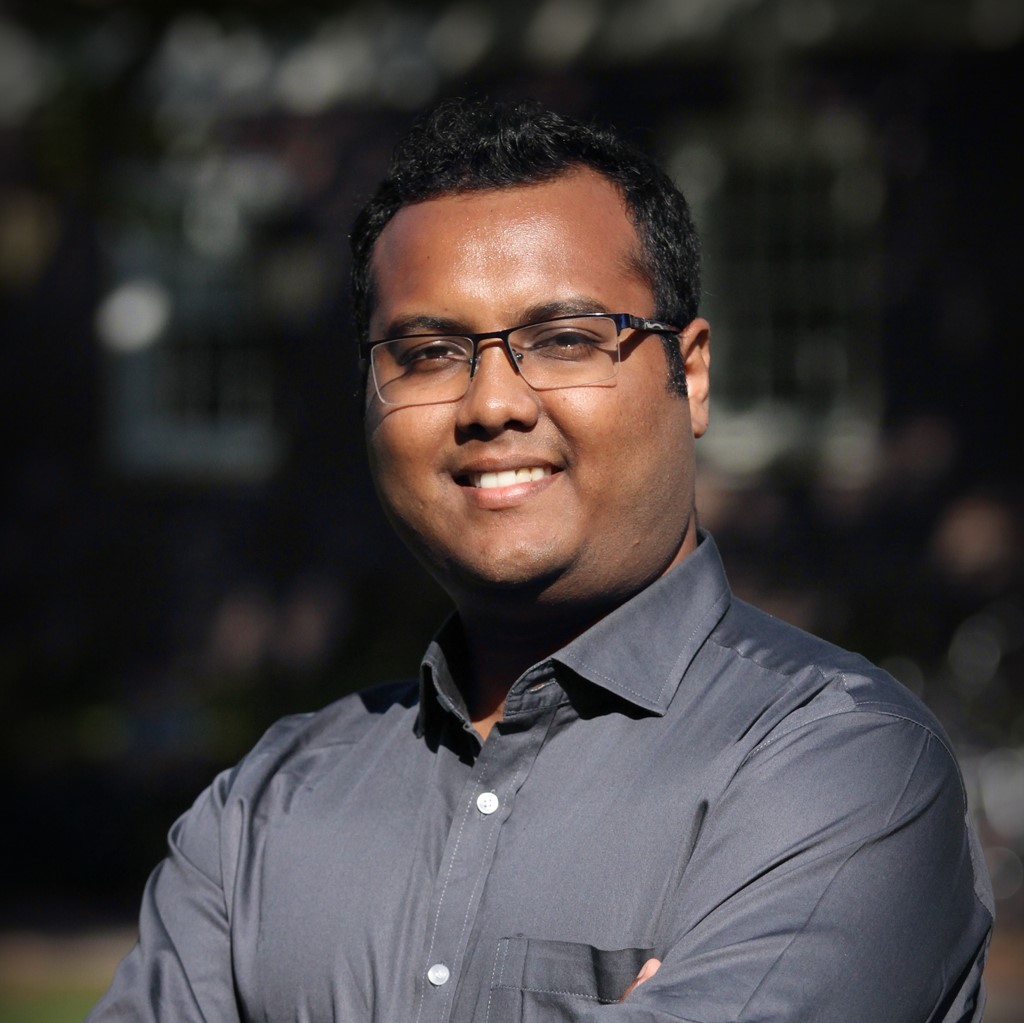}
\end{wrapfigure}
\textbf{A M Ishtiaque Mahbub} received his B.Sc. degree in mechanical engineering from Bangladesh
University of Engineering and Technology, Bangladesh in 2013 and a M.Sc. in computational
mechanics from University of Stuttgart, Germany in 2016. He is currently pursuing his Ph.D.
degree in mechanical engineering at the Information Decision Science Laboratory under the
supervision of Prof. Andreas A. Malikopoulos. His research interests includes, but are not limited
to, optimization and control with an emphasis on applications related to connected automated
vehicles, hybrid electric vehicles, and intelligent transportation systems. He has conducted
several internships at Robert Bosch LLC (USA), Robert Bosch GmbH (Germany), and Fraunhofer
IPA (Germany).

\par
 \begin{wrapfigure}{L}{0.17\textwidth}
\centering
\includegraphics[width=0.21\textwidth]{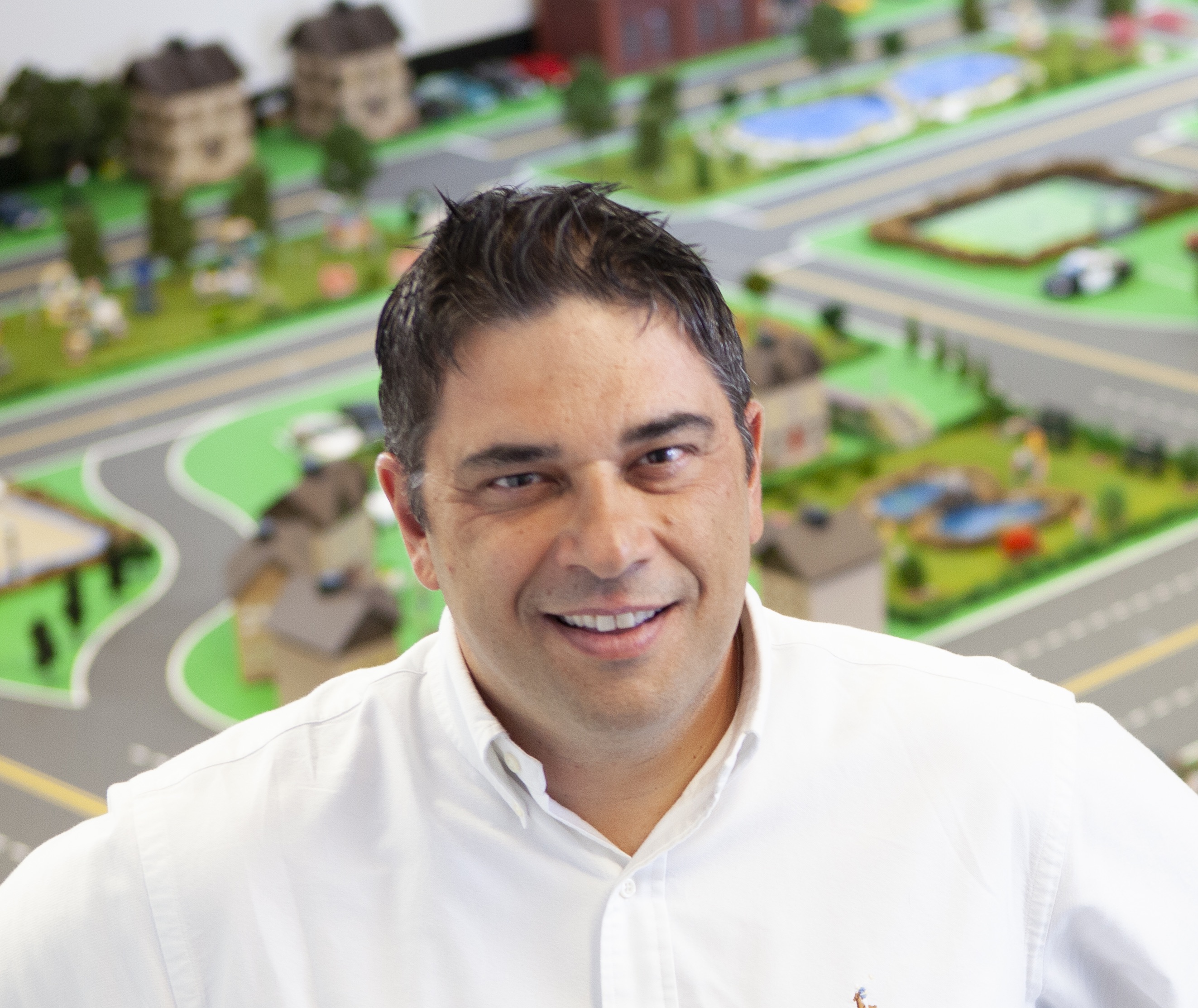}
\end{wrapfigure}
\textbf{Andreas A. Malikopoulos} received the Diploma in mechanical engineering from the National Technical University of Athens, Greece, in 2000. He received M.S. and Ph.D. degrees from the department of mechanical engineering at the University of Michigan, Ann Arbor, Michigan, USA, in 2004 and 2008, respectively. 
He is the Terri Connor Kelly and John Kelly Career Development Associate Professor in the Department of Mechanical Engineering at the University of Delaware (UD), the Director of the Information and Decision Science Laboratory, and the Director of the Sociotechnical Systems Center. Before he joined UD, he was the Deputy Director and the Lead of the Sustainable Mobility Theme of the Urban Dynamics Institute at Oak Ridge National Laboratory, and a Senior Researcher with General Motors Global Research \& Development. His research spans several fields, including analysis, optimization, and control of cyber-physical systems; decentralized systems; and stochastic scheduling and resource allocation problems. The emphasis is on applications related to sociotechnical systems, energy efficient mobility systems, and sustainable systems. He is currently an Associate Editor of the IEEE Transactions on Intelligent Vehicles and IEEE Transactions on Intelligent Transportation Systems. He is a member of SIAM, AAAS, a Senior member of IEEE, and a Fellow of the ASME.
\par
 \begin{wrapfigure}{L}{0.17\textwidth}
\centering
\includegraphics[width=0.21\textwidth]{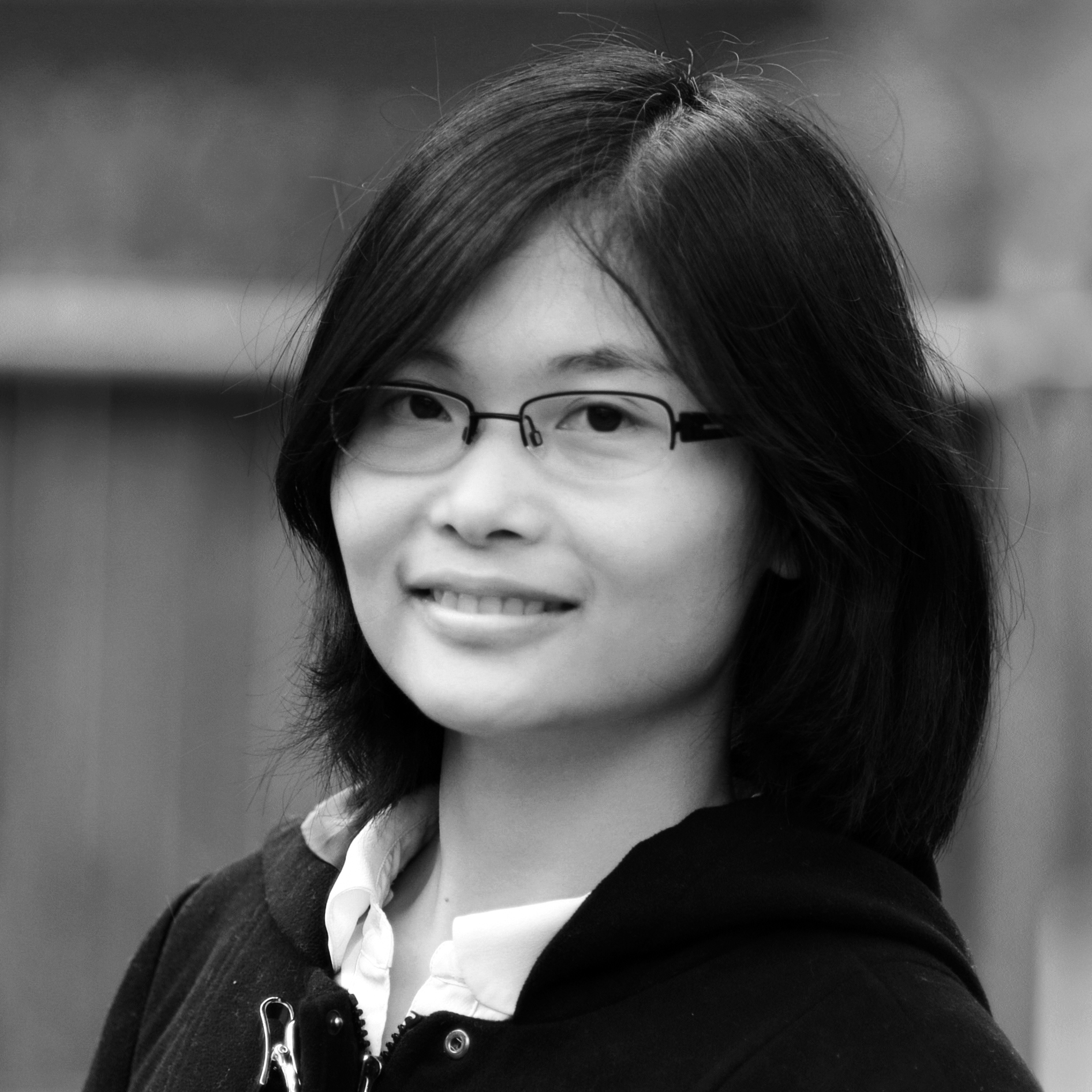}
\end{wrapfigure}
\textbf{Liuhui Zhao} received the B.S. degree
in Resources Science and Technology from Beijing
Normal University, Beijing, China, in 2009, the
M.S. degree from Department of Geography at the
University of Alabama in 2011, and the Ph.D. degree in Transportation Engineering from New Jersey
Institute of Technology in 2016. She is currently
a Postdoctoral Researcher in the Information and
Decision Science (IDS) Laboratory at the University
of Delaware leading research projects on emerging
transportation systems. She has participated in various research projects on connected automated vehicles, intelligent transportation systems, traffic and transit operations. Her research interests lie within
the areas of intelligent transportation systems, shared mobility, and connected
automated vehicles.
\end{document}